\documentclass[12pt]{amsart}
\usepackage{xspace,amssymb,amsfonts,epsfig,euscript,eufrak,xypic}
\author{Ben Webster}
\title[Stabilization phenomena in Kac-Moody algebras]
{Stabilization phenomena in Kac-Moody algebras and quiver varieties} 
\address{Department of Mathematics\\
         University of California, Berkeley\\
         970 Evans Hall\\
         Berkeley, CA 94720}
\subjclass[2000]{17B67}
\email{bwebste@math.berkeley.edu}
\urladdr{http://math.berkeley.edu/~bwebste}
\thanks{This material is based upon
  work supported under a National Science Foundation Graduate Research 
  Fellowship and partially supported by the RTG grant DMS-0354321.}
\begin{document}
\begin{abstract} Let $X$ be the Dynkin diagram of a symmetrizable
  Kac-Moody algebra, and $X_0$ a subgraph with all vertices of degree 1
  or 2. Using the crystal structure on the components of quiver
  varieties for $X$, we show that if we expand $X$ by extending $X_0$,
  the branching multiplicities and tensor product multiplicities stabilize, provided the weights involved 
  satisfy a condition which we call ``depth'' and are supported outside $X_0$.  This extends a theorem of Kleber and Viswanath.  

  Furthermore, we show that the weight multiplicities of such
  representations are polynomial in the length of $X_0$, generalizing
  the same result for $A_\ell$ by Benkart, {\em et al.}  
\end{abstract}
\maketitle

\newtheorem{thm}{Theorem}[section]
\newtheorem{prop}[thm]{Proposition}
\newtheorem{cor}[thm]{Corollary}
\newtheorem{lem}[thm]{Lemma}

\theoremstyle{remark}
\newtheorem*{defn}{Definition}
\newtheorem{remark}{Remark}

\newcommand{\nc}{\newcommand}
\newcommand{\renc}{\renewcommand}
  \nc{\kac}{\kappa^C}
  \nc{\Lam}[3]{\La^{#1}_{#2,#3}}
  \nc{\Lab}[2]{\La^{#1}_{#2}}
  \nc{\Lamvwy}{\Lam\Bv\Bw\By}
  \nc{\Labwv}{\Lab\Bw\Bv}
  \nc{\nak}[3]{\mathcal{N}(#1,#2,#3)}
  \nc{\hw}{highest weight\xspace} 
  \nc{\al}{\alpha}
  \nc{\be}{\beta}
  \nc{\bM}{\mathbf{m}}
  \nc{\bkh}{\backslash}
  \nc{\Bi}{\mathbf{i}}
  \nc{\Bj}{\mathbf{j}}
  \nc{\Bv}{\mathbf{v}}
  \nc{\Bw}{\mathbf{w}}
  \nc{\By}{\mathbf{y}}
  \nc{\Bz}{\mathbf{z}}
  \nc{\coker}{\mathrm{coker}\,}
  \nc{\C}{\mathbb{C}}
  \nc{\ch}{\mathrm{ch}}
  \nc{\de}{\delta}
  \nc{\ep}{\epsilon}
  \nc{\fr}[1]{\mathfrak{#1}}
  \nc{\GL}[2]{\mathrm{GL}_{#1} #2}
  \nc{\Hom}[3]{\mathrm{Hom}_{#3}(#1,#2)}
  \nc{\im}{\mathrm{im}\,}
  \nc{\La}{\Lambda}
  \nc{\la}{\lambda}
  \nc{\mult}{b^{\mu}_{\la_0}\!}
  \nc{\mc}[1]{\mathcal{#1}}
  \nc{\om}{\omega}
  \nc{\qvw}[1]{\La(#1 \Bv,\Bw)}
  \nc{\Rperp}{R^\vee(X_0)^{\perp}}
  \nc{\si}{\sigma}
  \nc{\SL}[1]{\mathrm{SL}_{#1}}
  \renc{\th}{\theta}
  \nc{\vp}{\varphi}
  \nc{\wt}{\mathrm{wt}}
  \nc{\Z}{\mathbb{Z}}
  \nc{\Znn}{\Z_{\geq 0}}
  \nc{\EV}{\EuScript{V}}
  \nc{\EE}{\EuScript{E}}
  \nc{\Spec}{\mathrm{Spec}}
  \nc{\tie}{\EuScript{T}}

\section{Introduction}
\label{sec:introduction}

\subsection{Background}
\label{sec:background}

If $X=(\EV,\EE)$ is the Dynkin diagram of a symmetrizable Kac-Moody
algebra, then we can construct a series of diagrams $X(m)$ by 
attaching an $A_m$ ``tail'' to a fixed vertex.  For large $m$, one can
consider weights supported on the ``head'' ($X$ and the adjacent portion
of the tail) and on the end of the tail, and vanishing on a large
section of tail in between.  We can extend such a weight to $X(m)$ for
all larger values of $m$ by simply allowing this gap to grow.  We might
hope that there would be some interesting asymptotic behavior as $m$
goes to infinity.

In the case where $X=A_m$, this construction was first studied by
Hanlon in \cite{Han85}, and R. Brylinski in \cite{Bry89}. They
discovered a number of interesting stabilization phenomena, most of
which seem specific to the $A_m$ case.  However, in \cite{KV04,Vis05},
Kleber and Viswanath show that one aspect does carry over:
\begin{thm}\label{KV-theorem}{\em (Kleber, Viswanath)}
  Let $\la,\mu,\nu$ be a deep triple of weights (as defined in Section
  \ref{sec:deep-weights}) on $X(m)$ supported on the
  head and end of tail as described above.  Then the tensor product 
  multiplicity $c^\la_{\mu,\nu}(m)$ for the Kac-Moody algebra $\fr g(m)=\fr
  g(X(m))$ stabilizes for $m$ sufficiently large.
\end{thm}

The proof is based on an analysis of the crystal graph of a
representation of $\fr g(m)$ via Littlemann paths and the following easy
consequence of the theory of crystal graphs due to Kashiwara and others
(see, for example, \cite{CP}) 
\begin{thm}\label{crystal-tens} Let $B_\mu$ be the crystal (graph)
  corresponding to the representation $V_\mu$. Then
  \begin{equation}
    c^\la_{\mu\nu}=\#\{x\in B_\mu|\mathrm{wt}(x)=\la-\nu, \tilde e_i^{\nu(\al_i^\vee)+1}x=0\},
  \end{equation}
  where $\tilde e_i:B_\mu\cup\{0\}\to B_\mu\cup\{0\}$ denotes the $i$th Kashiwara operator.
\end{thm}

Our aim in this paper is to generalize Theorem~\ref{KV-theorem} by 
using a variant of Kleber and Viswanath's techniques.  We will replace Littelmann 
paths with the crystal structure on the components of quiver 
varieties, described in \cite{Ks97,Sav04}.

This approach allows us to extend Theorem~\ref{KV-theorem} to a more general class
of arrangements. Our results imply the stabilization theorems of
\cite[Theorem 2.6]{KV04} and \cite[Theorem 2.5]{Vis05}, though some care
is required to see this in the case of \cite{KV04}: the restriction of
``extensibility'' there exactly implies that tensor product
multiplicities vanish in the non-deep case. Kleber and Viswanath refer to this as the
``number-of-boxes condition'' by analogy with $\mathfrak{gl}_n$, where it
implies that the number of boxes in the Young diagram adds under tensor
product.

\subsection{Summary of results}
\label{sec:summary-results}

Let $C$ be a generalized Cartan matrix, in the sense of Kac, and let
$C_0$ be a subset of the row/column set (identified by transpose) such
that each corresponding row or column contains at most 2 non-zero
off-diagonal entries, and all such entries are $-1$. For simplicity,
we also require that no minor consisting of rows from $C_0$ is the
Cartan matrix of $\tilde A_n$.

More geometrically, consider the directed graph $X$ with adjacency
matrix $A=2I-C$.  Throughout our paper, we will assume that $X$ is the
graph corresponding to a symmetrizable Cartan matrix. Then the subgraph
$X_0$ corresponding to $C_0$ must be a union of components
\begin{equation*}
  X_0=\bigsqcup_{i=1}^k X_i
\end{equation*}
with $X_i\cong A_{\ell_i}$ for some integers $\ell_1,\ldots,\ell_k$.
We call such a subgraph {\em elastic}.  The symbol $X_0$ will denote
an elastic subgroup throughout this paper.

\begin{figure}[htbp]
  \centering
    \centerline{\epsfig{figure=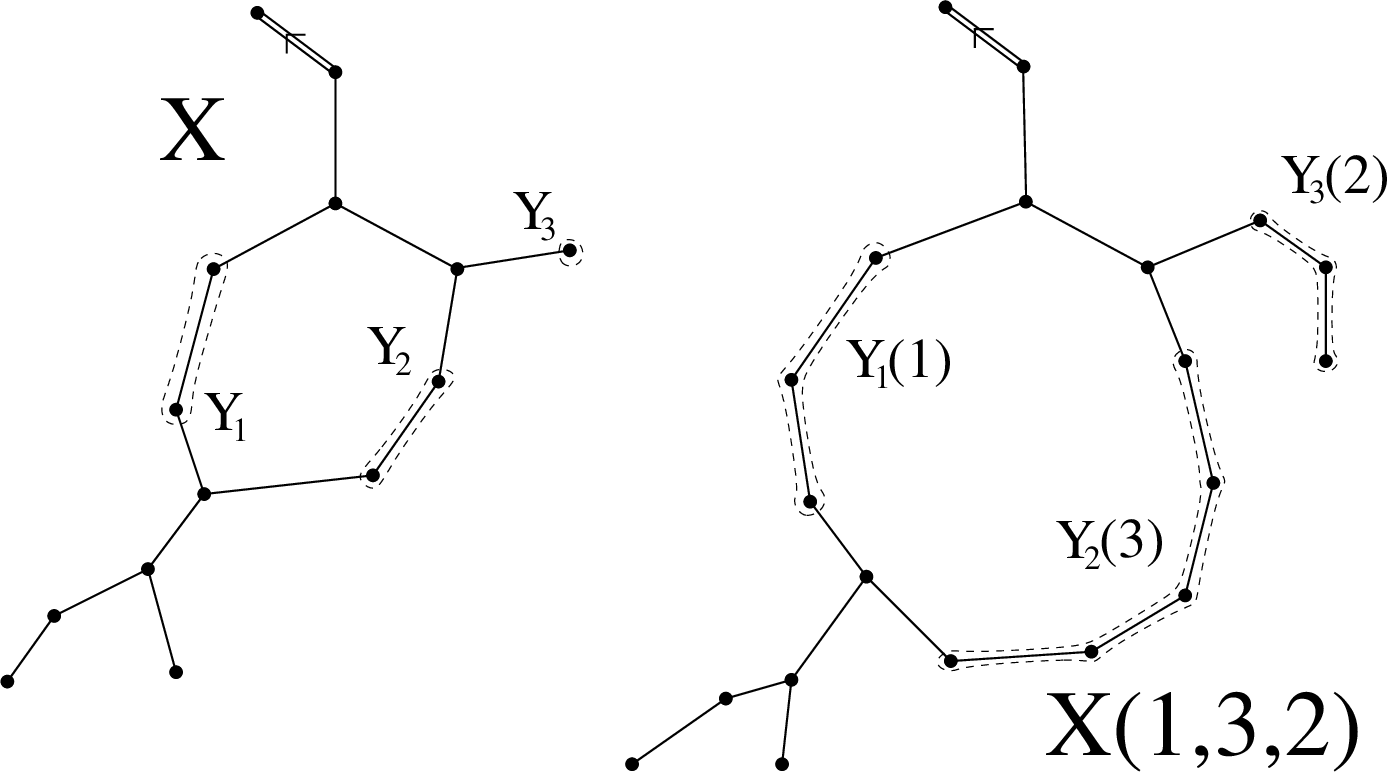, height=6cm}}
  \caption{The effect of expansion on a Dynkin diagram}
  \label{fig:1}
\end{figure}

We let $X(m_1,\ldots, m_k)=X(\bM)$ be $X$ with $X_0$ replaced by
\begin{equation*}
  X_0(\bM)=\bigsqcup_{i=1}^k X_i(m_i)\cong A_{m_i+\ell_i}.
\end{equation*}
That is, we let the selected string in our Dynkin diagram grow by $m_i$
edges.  This is illustrated in Figure~\ref{fig:1}.

Let $R^\vee (X_0)^\perp$ be the annihilator of the coroot lattice of
$X_0$, that is,  
\begin{equation*}
  R^\vee(X_0)^\perp=\{\de\in P(X)|\de(\al_i^\vee)=0, \,\forall i\in X_0\} 
\end{equation*}

If $\mu\in R^\vee(X_0)^\perp$, we will often say it ``vanishes on $X_0$.''

If we have a weight $\mu$ vanishing on $X_0$, we can view $\mu$ as a
weight of $X(m)$ for all $m$ (there are minor difficulties if the
Cartan matrix of $X$ is degenerate, but this will be dealt with in
later sections), since $X-X_0\cong X(\bM)-X_0(\bM)$.  Fix weights
$\mu,\nu,\la\in\Rperp$ and let $\th=\la-\nu$.  We let $V(\mu,\bM)$
denote the corresponding representation with \hw $\mu$ of the
Kac-Moody algebra $\fr g(\bM)$ with Dynkin diagram $X(\bM)$, and
$V(\mu,\bM)_\th$ its $\th$-weight space. Furthermore, we let $\fr
t(\bM)$ denote the Cartan subalgebra of $\fr g(\bM)$ and $\fr
g_0(\bM)$ the root subalgebra corresponding to the subdiagram $X_0$.
Note that $\fr {g_0+t}(\bM)$ is a finite dimensional reductive algebra
with the same Cartan subalgebra as $\fr g(\bM)$.  Thus $\mu$ also has
an associated irreducible highest weight representation of $\fr
{g_0+t}(\bM)$, which we denote $V_0(\mu,\bM)$.

Let us consider the asymptotic behavior of three of the most commonly
studied numerical invariants of representations of Kac-Moody algebras: 
\begin{itemize}
\item the branching multiplicities fro
m $\fr g(\bM)$ to $\fr g_0+\fr
  t(\bM)$: 
  \begin{equation*}
    b^{\mu}_{\th}(\bM)=\dim \Hom{V_0(\th,\bM)}{V(\mu,\bM)}{\fr{g_0+t}(\bM)}
  \end{equation*}
This is closely connected to the more often-studied branching rule for
$\fr g_0(\bM)$, but also includes the action of the full torus. 
This more detailed information will be important later. 
\item the tensor product multiplicities for $\fr g(\bM)$
  representations:
  \begin{equation*}
    c^\la_{\mu,\nu}(\bM)=\dim \Hom{V(\la,\bM)}{V(\mu,\bM)\otimes
    V(\nu,\bM)} {\fr g(\bM)}
  \end{equation*}
\item the weight multiplicities 
  \begin{equation*}
    w^\mu_\th(\bM)=\dim V(\mu,\bM)_\th.
  \end{equation*}
\end{itemize}

In the large $\bM$ limit, there is a subset of triples of weights
called {\em deep} which behave exactly as we might hope by
extrapolating from \cite{KV04} and the case of $A_\ell$:
\begin{thm} 
  \label{full-statement}
  For any elastic subgraph $X_0\subset X$ of the Dynkin graph of a
  symmetrizable Kac-Moody algebra, and for all weights $\la,\mu,\nu$,
  with $\mu+\nu-\la=\mu-\th$ deep of $j$-depth $\si_j$, and $m_i$
  sufficiently large for all $i\in \{1,\ldots, k\}$.
  \begin{enumerate}
  \item the branching multiplicities
  $b^{\mu}_{\th}(\bM)$ and tensor product multiplicities
  $c^{\la}_{\mu,\nu}(\bM)$ are independent of $\mathbf{m}$.
  \item the weight multiplicities $w_{\mu}^{\la}(\bM)$ are
  polynomials in the variables $m_i$ of multi-degree 
  $\leq (\si_1,\ldots,\si_k)$.
  \end{enumerate} 
\end{thm}

Unfortunately, these multiplicities remain very mysterious objects and
quite difficult to compute.  It would be particularly nice if a
more explicit combinatorial interpretation beyond Littelmann paths or
crystal graphs could be found, at least in some special cases. Even if
it were not computationally useful, it would represent a small window
into strange world of wild Kac-Moody algebras.

We begin with a brief introduction to deep weights in
Section~\ref{sec:deep-weights}, followed in
Section~\ref{sec:cryst-struct-quiv} by a review the relevant facts about
quiver varieties and their connections to representation theory which
will be necessary for the proof of this theorem.  Then we will prove
parts (1) and (2) of Theorem~\ref{full-statement} in
Section~\ref{sec:stab-quiv-vari} and
Section~\ref{sec:asympt-weight-mult} respectively.

\section*{Acknowledgments}
\label{sec:acknowledgements}

The author would like to thank Sankaran Viswanath, Joel Kamnitzer,
Alistair Savage and Nicolai Reshetikhin for their advice and help.  

We would also like to express our appreciation for the hospitality of
the Banff International Research Station and the efforts of the
organizers of the conference ``Representations of Kac-Moody Algebras and
Combinatorics'' held there in  March 2005, where this research
originated. 

\section{Deep weights}
\label{sec:deep-weights}

Let $P(X)$ be the weight lattice of $X$, and $R(X)\subset P(X)$ the root
lattice, and let $\al_i$ be the simple roots and $\al_i^\vee$
the simple coroots. Let $\Z^{\EV(X)}$ be the set of maps $\EV(X)\to\Z$.  There are 
natural maps
\begin{align*}
\om:P(X)&\to\Z^{\EV(X)}&\al:R(X) &\to \Z^{\EV(X)}\\
\la &\mapsto \{ i\mapsto \la(\al_i^\vee)\}&
\sum_{i\in \EV(X)}v_i\al_i &\mapsto \{i\mapsto v_i\}
\end{align*}
The map $\al$ is obviously an isomorphism, since the simple roots
$\al_i$ are linearly independent after choice of a realization (in the
sense of Kac).

The map $\om$ is surjective, but not in general injective.  Its kernel
is naturally identified with the nullspace of the Cartan matrix.  

Still, the properties of representations important to us depend only on
$\om(v)$.  If we have weights $\mu',\nu',\la'$ such that
\begin{align*}
  \Bw&=\om(\mu)=\om(\mu')\\
  \By&=\om(\nu)=\om(\nu')\\
  \Bv&=\al(\mu+\nu-\la)=\al(\mu'+\nu'-\la')
\end{align*}
then there is a unique isomorphism $B_\mu\cong B_{\mu'}$ between their
crystals (since $V(\mu)\cong V(\mu')$ as $\fr n^-$-modules), so
Theorem~\ref{crystal-tens} implies that
$c^\la_{\mu,\nu}=c^{\la'}_{\mu',\nu'}$.  Thus, we can denote this
multiplicity $c^{\Bv}_{\Bw,\By}$ without any loss of information.

Similarly, by the same crystal isomorphism, we can see that
$b^\mu_\th=b^{\mu'}_{\th'}$, so we can write
$b^{\mu}_{\th}=b^{\Bw}_{\Bv}$ where $\Bw=\om(\mu)$ and
$\Bv=\al(\mu-\th)$. 

Let $Z\subset X$ be any subgraph.  We call a function
$\Bv\in\Z^\EV(X)$ {\bf deep for $Z$} if $v_i$ is constant on each
connected component $Z_j$ of $Z$.  The value $\si_j$ of this constant
on $Z_j$ is called the {\bf ${Z_j}$-depth} of $\Bv$.  The {\bf depth}
of $\Bv$ is the maximum of the $Z_j$-depths.  In the case where
$Z=X_0$, we typically write $j$-depth for $X_j$-depth.

If $\la\in R(X)$, then $\la$ is called deep if $\al(\la)$ is deep.  A
triple of weights, $(\mu,\nu,\la)$ is called deep if $\mu+\nu-\la$ is
deep.  Depths and $j$-depths of weights are defined similarly.

\subsection{The case of Viswanath}
\label{sec:case-viswanath}

Depth seems like a strong condition, but in fact, non-deep
weights have a very hard time staying in the root lattice. In certain
cases ($A_n$, for example), the depth hypothesis in
Theorem~\ref{full-statement} is entirely redundant. Assume $X$ and
$X_0$ are connected, and $X\bkh X_0$ is disconnected, so $X_0$ is
connecting two separate Dynkin diagrams.  This situation has been
studied by Viswanath in \cite{Vis05}.  If these diagrams form an {\bf
extensible pair} as defined in \cite{Vis05} then for large $m$, the diagram $X(m)$ has
non-degenerate Cartan matrix, and the fundamental weights $\om_j$ make
sense.  Thus, we have an inclusion
$i_m:R^\vee(X_0)^{\perp}\hookrightarrow R(X_0(m))^{\perp}$ such that
$i_m(\om_j)=\om_{j'}$.

\begin{prop}\label{sec:pop-out}
  {\em (Viswanath)} If $\la\in R(X)$ is not deep, then
  $i_m(\la)\notin R(X(m))$ for all $m$ sufficiently large. Thus,
  $c^\la_{\mu,\nu}(m)$ stabilizes for all weights.
\end{prop}

While there are many examples of extensible pairs, there are a lot of
cases that are not extensible and Proposition~\ref{sec:pop-out} fails.
For example, while $X=A_n$ with $X_0$ connected is extensible in this
sense, types $B_n, C_n, D_n$ are not, and in these cases
counterexamples to Proposition~\ref{sec:pop-out} are easily
located. Interestingly, numerical evidence suggests that even in these
cases, tensor product multiplicities may stabilize, but beyond
computer computations and the finite and affine cases (as described in
\cite{BBL90,BKLS99}), little is known about these cases.

\begin{figure}
  \centering
    \centerline{\epsfig{figure=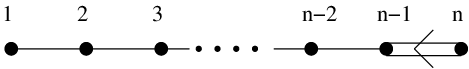, height=1.2cm}}
  \caption{The Dynkin diagram $B_n$}
  \label{fig:Bn}
\end{figure}

For example, let $X=B_n$, with the numbering given in
Figure~\ref{fig:Bn}. 

For all $n$,$V(\om_n)$ is the spin representation, and we have the
identity $$2\om_n=\sum_{i=1}^n i\al_i$$ so the labeling of
Figure~\ref{fig:Bn} also represents the values of $\al(2\om_n)$. It
follows that $2\om_n\in R(B_n)$ for all $n$, and is obviously 
not deep for any $X_0\subset X$    

Since $V(\om_n)\cong V(\om_n)^*$, for all $n$, $c^{0}_{\om_n,\om_n}=1$,
even though this triple is not deep.  Thus, in this case, and many
others, we cannot ignore non-deep weights.  

\section{Crystal structures on quiver varieties}
\label{sec:cryst-struct-quiv}

\subsection{Quivers and their representations}
\label{sec:quiv-their-repr}

Our primary tool will be the geometry of quiver varieties and its 
connections with representation theory, originally due to Lusztig 
and Nakajima.  See the papers \cite{Lus91,Nak94,Nak98} for more 
details on the construction of quiver varieties and their applications.

A {\bf quiver} is a directed graph $X=(\EV,\EE)$.  A {\bf framed
  representation} of a quiver is a choice of:
\begin{itemize}
\item two vector spaces, $V_i$ and $W_i$ for each vertex $i\in\EV$,
\item maps $x_e:V_{\al(e)}\to V_{\om(e)}$ and $x_{\bar e}:V_{\om(e)}\to 
V_{\al(e)}$ for each $e\in\EE$ (here we use $\al(e)$ and
$\om(e)$ to denote the
tail and head of the edge $e$).
\item maps $t_i:V_i\to W_i$ for each vertex $i\in\EV$.
\end{itemize}
It is often convenient to package together the vector spaces 
\begin{equation*}
  V=\bigoplus_{i\in \EV(X)} V_i \qquad W=\bigoplus_{i\in \EV(X)} W_i
\end{equation*}
and view them as a representation $V\oplus W$ of the framed
path algebra $\Pi(X)$.  This algebra is generated by elements
$x_P$ for each path $P$ in $X$, and $j_i$ for $i\in\EV(X)$, with
relations:
\begin{align*}
  x_Px_{P'}&=x_{P*P'},& j_ij_k&=0 \text{ for all } i,k\in \EV\\
  x_Pj_i&=0, & j_ix_P&=0 \text{ iff } i\neq \om(P).\\
\end{align*}
The vector space $V\oplus W$ is naturally a representation of this algebra with $x_P$
given by composition of $x_e$ for the edges $e\in P$, and the action of
$j_i$ given by the map $t_i$.

Let $\mc I$ be the ideal generated by all paths of length greater than
0, and let $\ep_i(x,t)=\dim (V/\mc I V)_i$ be the dimension of the
quotient of the vector space $V_i$ by the image of all the paths
ending at $i$.

We say a quiver representation
\begin{itemize}
\item is {\bf nilpotent} if for some integer $N$, $x_P=0$ for all paths $P$
  with length $\geq N$.
\item is {\bf stable} if there is no subrepresentation contained in
  $V$.
\item satisfies {\bf property MM} if 
  \begin{equation*}
    \sum_{e\in \EE(X)} x_ex_{\bar e}-x_{\bar e}x_e=0.
  \end{equation*}
\end{itemize}

\begin{defn}
  A {\bf Nakajima quiver representation} (or NQR for short)
is a stable, nilpotent framed quiver
representation which satisfies property MM, along with a system of isomorphisms
$\vp_i:W_i\to\C^{w_i}$.

Let $\nak XVW$ be the variety of NQR's for $X$ on the vector space
$V\oplus W$.
\end{defn}

We let $\La(X,V,W)$ be the quotient of the variety of NQR's
on $V\oplus W$ under the obvious 
action of 
\begin{equation*}
G=\prod_{i\in\EV}\mathrm{GL}(V_i)\times \mathrm{GL}(W_i).    
\end{equation*}
The condition of stability implies that $G$ acts freely on $\La(X,V,W)$,
so this is a geometric quotient, and using geometric invariant theory,
we can give it the structure of a projective variety over $\C$.

This is a Lagrangian subvariety inside the
hyperk\"ahler quiver variety defined by Nakajima in \cite{Nak94}.  In
fact it, it is the ``core,'' in the sense of Proudfoot \cite[\S 2.2]{Pro04}.

Up to natural isomorphism, this variety only depends on the dimensions 
\begin{equation*}
v_i= \dim V_i \qquad 
w_i= \dim W_i  
\end{equation*}
Thus, we will usually denote this variety by
$\qvw{X,}$ without choosing fixed vector spaces of the given
dimensions.

\subsection{Quiver varieties as moduli spaces}
\label{sec:quiver-varieties-as}

Let $\La(X)$ be the union over all $\Bv$ and $\Bw$ of
$\La(X,\Bv,\Bw)$.  Then $\La(X)$ also has a moduli theoretic
interpretation:

Let $S$ be a scheme over $\C$.  Then our notion of a NQR of vector
spaces can be extended to an NQR of algebraic vector bundles
(i.e. locally free sheaves) over $S$.  Following the usual terminology
of algebraic geometry, we refer to this as a family of NQR's over $S$.
The case of vector spaces is simply when $S=\Spec\, \C$.

Over $\nak XVW$, there is a natural, $G$-equivariant family
$(\mc{\tilde{V} \oplus \tilde{W}},\tilde{x},\tilde{t})$ of NQR's, which
one may call {\bf the tautological NQR}.  Since
$G$ acts freely on $\nak XVW$, this is the pullback of an NQR
$(\mc{V\oplus W},x,t)$ on $\qvw{X,}$, and taking union, on $\La(X)$.  We
call this {\bf the universal NQR}.

Consider the cofunctor $\mc N:\mathrm{Sch}/\C\to \mathrm{Set}$ sending a scheme
over $\C$ to the set of isomorphism classes of families of NQR's.  As
usual, pullback of the universal bundle by $f$, gives a natural
transformation of functors $U:\mathrm{Mor}(-,\La(X))\to\mc N$.  
\begin{prop}
  The natural transformation $U$ is an isomorphism.  That is, the scheme
  $\La(X)$ represents $\mc N$.   
\end{prop}
\begin{proof}
  Since we can glue morphisms, we can reduce to the case where $S$ is
  affine and our NQR is on a trivial vector bundle $L_V\oplus L_W$.
  $L_W$ is already trivialized by the definition of an NQR.  Thus, we
  need only pick a basis of nonvanishing sections of $L_V$.  Taking matrix
  coefficients, this gives a map from $S$ to $\nak
  X{\oplus\C^{v_i}}{\oplus\C^{w_i}}$ and thus by projection to
  $\La(X,\Bv,\Bw)$.  Obviously, the tautological bundle on $\nak
  X{\oplus\C^{v_i}}{\oplus\C^{w_i}}$ pulls back to $L_V\oplus L_W$.  By
  transitivity of pullback, the universal bundle does as well.

  Of course, changing the basis results in the same map to $\La(X,\Bv,\Bw)$
  via a different lift to the space of NQR's.  Thus the map we have
  constructed is unique. 
\end{proof}

\subsection{Connections with Kac-Moody algebras}
\label{sec:connections-with-kac}

Now assume that $X$ is simply laced (the corresponding Cartan matrix entries
are all in $\{2,0,-1\}$), and let $\fr g$ be
the Kac-Moody algebra associated to the graph $X$.  The geometry of
the quiver variety $\La$ interacts with the representation theory of
$\fr g$ in a number of interesting ways. For us, the important one
will be the quiver model for the crystal $B_\mu$ described in
\cite{Sav04}.

Choose an orientation on $X$ with no oriented loops (for example, pick
an ordering on the vertices, and direct the edges in ascending order).
Let $\mu$ be a weight for $X$ with $\om(\mu)=\Bw$, 
  (i.e. $\mu(\al_i^\vee)=w_i$).  Then we have the following remarkable
  result from \cite{Saito02}
\begin{thm}\label{crys-state} {\em (Saito)}
  There is a natural $\fr g$-crystal structure 
  $C_{\Bw}$ on the components of 
  \begin{equation*}
    \La(X,\Bw)=\bigsqcup_{\Bv\in \Z_{\geq 0}^{\EV(X)}}\qvw{X,}
  \end{equation*}
  such that  
  \begin{enumerate}
  \item $C_{\Bw}\cong B_\mu$.
  \item The weight map $\mathrm{wt}:C_{\Bw}\to P(X)$ is given by 
  \begin{equation*}
    \wt(c)=\mu-\al^{-1}(\Bv)
  \end{equation*}
  where $c\subset \qvw{X,}$.
  \item The set of representations $\{(V,W,x,t)\vert\ep_i(x,t)\leq
  n\}$ is dense in a component $c$ if and only if $\tilde
  e_i^{n+1}c=0$, where, as before, $\tilde e_i$ is a Kashiwara
  operator.
\end{enumerate}
\end{thm}

Thus, the statistics of elements of the crystal which are important to
us in determining tensor product multiplicities and branching rules have
clear interpretations in terms of the associated quivers.

Consider the open subsets of $\La(X)$ defined by
\begin{align*}
  \Lab\mu\th(X,X_0)&=\{(V,W,x,t)\in\La(X,\om(\mu),\al(\mu-\th))|\ep_i(x,t)=0
  \text{ for all } i\in X_0\}\\
\Lam\la\mu\nu (X)&=\{(V,W,x,t)\in\La(X,\om(\mu),\al(\mu+\nu-\la))|\ep(x,t)\leq\om(\nu)\}.   
\end{align*}
These can be interpreted as moduli spaces of NQR's with fixed
representation theoretic properties.

As an immediate consequence of Theorem~\ref{crys-state}, we see the following proposition.
\begin{prop}\label{bran-ten-comp} For all weights $\la,\mu,\nu,\th$,  
  \begin{itemize} 
  \item The branching multiplicity $b^\mu_\th(X,X_0)$ is the number of
  components of $\Lab\mu\th(X,X_0)$.

  \item The tensor product multiplicity $c^{\la}_{\mu,\nu}$ is the number of
  components of $\Lam\la\mu\nu (X)$.
  \end{itemize}
\end{prop}

As was discussed in Section~\ref{sec:deep-weights}, we will also use the
notation $\Lab{\Bw}{\Bv}$ and $\Lam\Bv\Bw\By$, where $\Bw=\om(\mu)$,
$\By=\om(\nu)$ and $\Bv=\al(\mu+\nu-\la)=\al(\mu-\th)$, interchangeably with
$\Lab\mu\th$ and $\Lam\la\mu\nu$.

If $X$ is not simply laced, but symmetrizable, then we let
$\qvw{X,}$ be the appropriate closed subvariety of
$\La(\tilde{X},\tilde{\Bv},\tilde{\Bw})$, where $\tilde X$ is a simply
laced diagram with an automorphism $a$ such that $\tilde X/a\cong X$
(see \cite{Sav04} for the details of this construction for arbitrary
symmetrizable algebras).  This corresponds to the well-known embedding
of $B_\mu$ into the corresponding crystal graph for $\tilde X$.  

\section{Stabilization of quiver varieties}
\label{sec:stab-quiv-vari}

\subsection{Notation}
\label{sec:notation}

Throughout this section, we fix weights  $\mu,\nu,\la$, with
$\mu+\nu-\la\in R(X)$. We let $\Bw=\om(\mu)$, $\By=\om(\nu)$,
$\Bv=\al(\mu+\nu-\la)$, and $\la=\nu+\th$.  

Let $X$ and $X_0$ be as in the introduction, and assume that $\Bw$ and
$\By$ vanish on $X_0$.  From the definition, it follows that
$\Lamvwy(X)\subset \Labwv(X,X_0)$.

Fix an integer $s$, and let $X_0^s$ be the set of vertices in $X_0$
which are more than $s$ edges from the border of $X_0$, by which we
mean the set of vertices in $X_0$ which are adjacent to edges in
$X\bkh X_0$.  For the remainder of this section, we assume that
$X_0^{s+1}$ is non-empty.

Fix a vertex $i\in \EV(X_0^{s+1})$. By renumbering, we may assume that
$i\in \EV(X_1)$. Let $X'$ be a graph with vertex set $\left(\EV(X)\bkh
\{i\}\right)\cup\{i',i''\}$ and edge set $\EE(X)\cup \{d\}$, with one
edge connected to $i$ attaching to $i'$ and the other to $i''$ and $d$
connecting $i'$ with $i''$.

We
let $X_0'$ be the obvious subgraph of $X'$ which is also a disjoint
union of type $A$ diagrams.  This addition of an edge is shown in
Figure~\ref{fig:2}. 

Now assume $\Bv$ is deep for $X_0^s$. It extends to a deep function
$\Bv'$ on $X'$ in an obvious, unique way; outside of $X_0'^s$, it is
defined by the graph isomorphism $X-X_0^s\cong X'-X'^s_0$, and on
$X_0'^s$ it simply takes the value of the $1$-depth $\sigma_1$ on $i'$
and $i''$, as is necessary to preserve deepness. We can extend $\Bw$ and
$\By$ to $\Bw'$ and $\By'$ by the same process, since they are also deep
for $X_0^{s+1}$.

\subsection{Stabilization}
\label{sec:stab-quiv-vari-1}

As before, $X$ is the Dynkin graph of a symmetrizable Kac-Moody
algebra, and $X_0\subset X$ is an elastic subgraph.
\begin{thm}\label{quiv-iso}
  Assume $v_i\leq s$ for all $i\in
  \EV(X_0)$.  If $\Bv$ is deep for $X_0^{s}$, there are canonical
  isomorphisms, preserving $\ep_j$ for all $j\in \EV(X)$, such that 
  \begin{equation*}
    \Labwv(X,X_0)\cong\Lab{\Bw'}{\Bv'}(X',X_0')
  \end{equation*}
  These in turn induce isomorphisms
  \begin{equation*}
        \Lamvwy(X)\cong\Lam{\Bv'}{\Bw'}{\By'}(X')
  \end{equation*}
  If $\Bv$ is not deep for $X_0^s$, then $\Labwv(X,X_0)$, and hence
  $\Lamvwy(X)$, is empty. 
\end{thm}

\begin{cor} \label{vanishing-cor}
  If $b^\mu_\th(X, X_0)\neq 0$ (resp. $c^\la_{\mu,\nu}(X)\neq 0$) and 
  $\al(\mu-\th)\leq s$ on $X_0$, then $\mu-\th$ is deep for $X_0^{s+1}$.
\end{cor}

\begin{cor} \label{branching-cor}
   $ \displaystyle{b_{\Bv}^{\Bw}( X, X_0)=b_{\Bv'}^{\Bw'}( X', X_0').}$
\end{cor} 
\begin{cor} \label{tensor-cor}
     $ \displaystyle{ c^\Bv_{\Bw,\By}( X)=c^{\Bv'}_{\Bw',\By'}( X').}$
\end{cor}

Applying these corollaries inductively, we can lengthen the strings in
$X_0$ arbitrarily without affecting branching or tensor product
multiplicities once $X_j(\bM)$ has at least $2\si_j+1$ edges. Thus,
Corollaries~\ref{branching-cor}~and~\ref{tensor-cor} imply part (1) of
Theorem~\ref{full-statement}. 

These corollaries seem to be stronger than the statement of
Theorem~\ref{full-statement}, since in this case $\la$ is not assumed to
vanish on $X_0$ but on $X_0^{s}$.  But in fact, as
Corollary~\ref{vanishing-cor} shows, $\la$ must vanish on $X_0^{s}$, and
we lose nothing in the large $m$ limit by shrinking $X_0$ by a fixed
number of vertices, so Theorem~\ref{full-statement} is simply the
asymptotic statement of
Corollaries~\ref{branching-cor}~and~\ref{tensor-cor}.

\begin{proof}[Proof of Theorem~\ref{quiv-iso}]
The idea of the proof is as follows: We show that on $X_0^{s}$, any
stable quiver representation must have a very simple form, so there is a
canonical way to extend it to an extra edge, or to remove an edge. 

First, we assume that $X$ is simply-laced.

We must define some notation. Let $d_{i}$ and
$d_{i+1}$ be the edges adjacent to $i$, and let $i+1$ be the vertex at
the opposite end of $d_{i+1}$ (similarly for $i-1$, $i\pm n$, etc.).  In
$X'_0$, we let $i'$ be adjacent to $i-1$ and $i''$ to $i+1$.

For a quiver representation $(x,t)$ on $X$, we denote by $x_{d_i}$ the 
quiver map $V_i\to
V_{i-1}$, and by $y_{d_i}=x_{\bar{d}_i}$ the map from $V_{i-1}\to V_i$. 
In the notation of Section~\ref{sec:cryst-struct-quiv}, this corresponds
to choosing an orientation of $X_0$ with directed edges pointing in the
decreasing direction.  

We use $x$ and $y$ instead of just $x$ in order to simplify notation in
the proof. When the point of application is obvious, we will drop the
subscript from the maps $x_{d_i}$ and $y_{d_i}$.  The notation defined
above is shown in Figure~\ref{fig:2}.

\begin{figure}
  \centering
    \centerline{\epsfig{figure=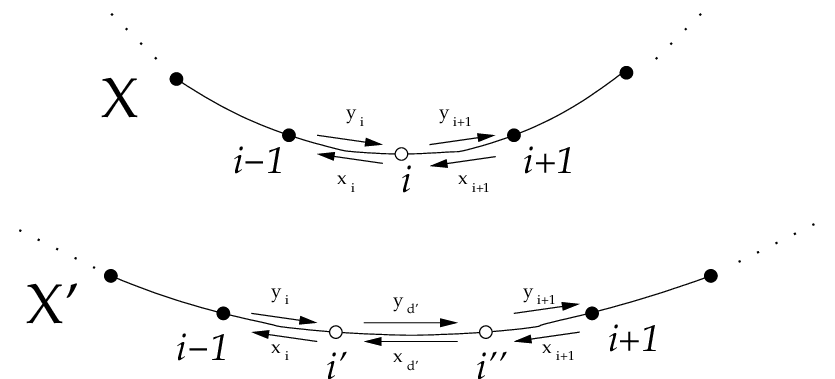, height=5cm}}
  \caption{The addition of an edge}
  \label{fig:2}
\end{figure}

If $(x,y)\in\qvw{X,}$, it must satisfy the moment map condition, which
reduces in the bivalent case to $y_{d_i}x_{d_i}=
x_{d_{i+1}}y_{d_{i+1}}$.  Now, apply this relation inductively to obtain
the equalities
\begin{equation}\label{moment-map}
  y^sx^s=(y_{d_i}x_{d_i})^s=
  x^sy^s 
\end{equation}

By assumption, $y_{d_i}x_{d_i}:V_i\to V_i$ is nilpotent. Since $v_i \leq s$,
$(y_{d_i}x_{d_i})^s=0$, and by \eqref{moment-map}, $x^s y^s=y^s x^s=0$ as well.

\renewcommand{\labelitemi}{$(*)$}

On the other hand, we claim: 
\begin{itemize}
\item $(\ker x^s)_j\cap(\ker y^s)_j=0$. 
\end{itemize}
If $z\in (\ker y^s)_j\cap(\ker x^s)_j$ and $z\neq 0$, then for some $\ell<s$,
$y^\ell z\neq 0$, but $y^{\ell+1}z=0$.  By assumption,  
\begin{equation*}
  x^s y^{\ell}z= y^{\ell} x^sz=0,
\end{equation*}
so there exists an integer $k$, such that $z'=x^k y^\ell z\neq 0$, but
$xz'=0$ and 
\begin{equation*}
  yz'=x^k(y^{\ell+1}z)=0.
\end{equation*}
This contradicts stability, so the claim $(*)$ is true.

\renc{\labelitemi}{$(**)$}

We claim also that for all $j\in \EV(X_0)^s$: 
\begin{itemize}
\item $V_j=(\im x^s)_j+(\im y^s)_j$.  
\end{itemize}
We prove this by induction on $s$.  For $s=1$, this is simply the
condition $\ep_j((x,y),t)=0$.   

Assuming the claim for $k-1$, and applying the relation 
\begin{equation*}
  (\im x^{\ell_1}y^{\ell_2})_j=(\im x^{\ell_1+1}y^{\ell_2})_j+
  (\im x^{\ell_1}y^{\ell_2+1})_j
\end{equation*}
several times, we see that 
\begin{align*}
  V_j&=(\im x^{k-1})_j+(\im y^{k-1})_j\\
     &=(\im x^k)_j+(\im y^k)_j+(\im x^{k-1}y)_j+(\im xy^{k-1})_j\\
     &=(\im x^k)_j+(\im y^k)_j+(\im x^{k-1}y^2)_j+(\im x^2y^{k-1})_j\\
     &\hspace{1.24in}\vdots\\
     &=(\im x^k)_j+(\im y^k)_j+(\im x^{k-1}y^k)_j+(\im x^ky^{k-1})_j\\
     &=(\im x^k)_j+(\im y^k)_j
\end{align*}

This proves the claim $(**)$.

Since $x^sy^s=0$, $(\im x^s)_j\subseteq (\ker y^s)_j$ and $(\im y^s)_j\subseteq(\ker
x^s)_j$. Applying $(*)$ and $(**)$, we see that $V_i=(\im x^s)_i\oplus(\im
y^s)_i=(\ker y^s)_i\oplus(\ker x^s)_i$.  

In particular, if $i-1$ is in $\EV(X_0)^s$, $x_{d_i}$ is an isomorphism between
$(\im x^s)_i\subset V_i$ and $(\im x^s)_{i-1}\subset V_{i-1}$ and $y_{d_i}$ an
isomorphism between $(\im y^s)_{i-1}\subset V_{i-1}$ and $(\im y^s)_i\subset V_i$,
and $v_i=v_{i-1}$. Applying this to all vertices in $X_0^s$, we see that
if there is a stable quiver in $\qvw{X,}$, then $v_i$ is constant on
connected components of this set, implying the second part of our
theorem. 

This shows that on $X_0^{s}$, our quiver representation has a very
simple structure.  It naturally decomposes into two subrepresentations,
one on which $x$ is an isomorphism, and one on which $y$ is an isomorphism, and the moment
map condition guarantees that our representation will resemble
Figure~\ref{fig:3}.  This figure should be read such that every solid dot
represents a basis vector of $V_i$, with leftward arrows representing
the action of $x$ and the rightward $y$.  Compare these to the case of
$A_\ell$, covered by Frenkel and Savage in \cite{FS03, Sav03} 
\begin{figure}
  \centering
        \centerline{\epsfig{figure=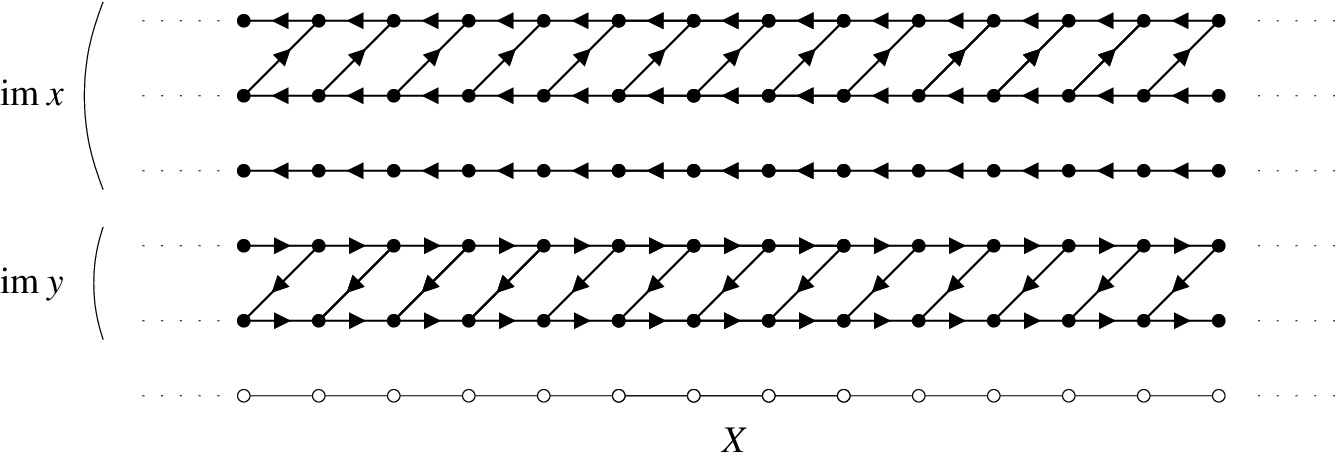, height=4.5cm}}
  \caption{The structure of a quiver representation on $X_0^s$.}
  \label{fig:3}
\end{figure}

From this picture, it is clear how to extend the pattern to the
additional edge $d$. To be precise, we let $\be':V_{i'}\to V_i$ be an
isomorphism, and similarly for $\be''$.  The map $x_d: V_{i''}\to
V_{i'}$ is defined by $(\be')^{-1}\be''$ on $(\be'')^{-1}(\im x^s)_i$
and $(\be')^{-1}x_{i+1}y_{i+1}\be''$ on $(\be'')^{-1}(\im y^s)_i$, and
the same for $y$, {\em mutatis mutandis}.

We had to make an arbitrary choice of $\be'$ and $\be''$, but we
obtain an isomorphic NQR for any choice.  In terms of our moduli
theoretic interpretation, this says that we have a natural NQR of vector
bundles for $X'$ on $\Labwv(X,X_0)$, and thus a natural map  
$\vp_{X, X'}:\Labwv(X,X_0)\to\Lab{\Bw'}{\Bv'}(X',X_0')$.

This map is in fact an isomorphism, since it has a clear inverse; given
a quiver of the same form on $X'$, the maps $x$ and $y$ are isomorphisms
on $(\im x^s)_i$ and $(\im y^s)_i$, so we can just as easily perform the
inverse operation of contracting an edge.

Note that we should not simply compose the quiver maps, as this will
violate the the moment map condition, but instead must use $x$ and $y$
to identify the vector spaces on adjacent vertices.

Thus, we have a canonical bijection
$\Labwv(X,X_0)\to\Lab{\Bw'}{\Bv'}(X',X_0')$.

Any reader with concerns about whether this bijection is in fact an
algebraic isomorphism, need only note that the process of adding or
contracting an edge can be performed in families, so it gives an
isomorphism between the functors representing $\Labwv(X,X_0)$ and
$\Lab{\Bw'}{\Bv'}(X',X_0')$.  Since the scheme representing a functor is
unique up to isomorphism, these varieties are canonically isomorphic. 

This isomorphism preserves $\ep_j$ for all $j\in \EV(X)$, since if $j$ is 
not in one of the $A_n$ strings, then the local structure of the quiver 
representation is unchanged.  Thus it also induces an isomorphism 
\begin{equation*}
  \Lamvwy(X)\cong\Lam{\Bv'}{\Bw'}{\By'}(X').
\end{equation*}

If $X$ is not simply laced, we simultaneously apply the result to all 
vertices covering $i$ in the simply laced diagram $\tilde X$, and 
deduce the same result. 
\end{proof}

\section{The polynomial behavior of weight multiplicities}
\label{sec:asympt-weight-mult} 
As another consequence of our quiver theoretic picture of stabilization,
we consider the behavior of weight multiplicities in the large $m$ limit
of our sequence of diagrams.

As with tensor product and branching rules, the dependence of weight
multiplicities on the weights involved can be rewritten in terms of
functions $\EV(X)\to\Z$.  We let $w^\mu_\th=w^{\Bw}_{\Bv}$, using the
notation of Section~\ref{sec:deep-weights}. 

The behavior of weight multiplicities is considerably more subtle than
that of branching rules.  For example, the branching rule from $\mathfrak{sl}_m$
to the Levi component of a parabolic is simple, well-known, and
stabilizes as $m\to\infty$, as is clear from their description in terms
of tableaux, or partial Gelfand-Tsetlin patterns. On the other hand, the
Weyl character formula is computationally quite unwieldy. However, the
large $m$ asymptotics of weight multiplicities for $\mathfrak{sl}_n$ are
reasonably well understood.  

\begin{thm}\label{benkart-et-al}
  \emph{(\cite{BKLS99})} Let $X\cong A_\ell$ and assume $\Bw$ vanishes
  on some connected subdiagram $X_0$ and $\Bv$ is of depth $s$ on $X_0$.
  Then for $m$ sufficiently large, $w^{\Bw}_\Bv(m)$ is equal to a
  polynomial in $m$ of degree $\leq s$.   
\end{thm}

Now, using the results of the previous section, we can reduce the
general case to the above result for $A_\ell$.  

As in the introduction, we can decompose $X_0(\bM)$ into its connected 
components 
\begin{equation*}
  X_0(\bM)=\bigsqcup_{i=1}^kX_i(m_i)\cong A_{\ell_i+m_i}
\end{equation*}
and let $\si_j$ be the $j$-depth of $\Bv$.
\begin{thm}\label{polyn-behav-weight}
  Let $X,X_0,\Bw,\Bv$ be as in Theorem~\ref{quiv-iso}.
  Then, for $\bM$ sufficiently large, the weight multiplicity $w^\Bw_\Bv(\bM)$
  is a polynomial in the variables $m_1,\ldots, m_k$ of multidegree less
  than $\boldsymbol{\si}=(\si_1,\ldots,\si_k)$. 
\end{thm}
This immediately implies Theorem~\ref{full-statement}, part (2).

\begin{lem}\label{weight-decomp}
   $ \displaystyle{w^\mu_\th( X)=\sum_{\xi\in R^+(X_0)} 
      b^{\mu}_{\th+\xi}(X,X_0) \cdot w^{\th+\xi}_{\th}( X_0).}$
\end{lem}

\begin{proof} 
  Let $\fr g$ and $\fr g_0$ be the Kac-Moody algebras corresponding to
  $X$ and $X_0$ respectively, and let $\fr t$ be a Cartan subalgebra
  of $\fr g$.  Any integral weight $\mu$ of $X$ can be considered
  as a integral weight on the reductive Lie algebra $\fr g_0+\fr t$.
  Let $Z_{\mu}$ be the corresponding finite dimensional highest weight
  representation. As with any reductive group,
  \begin{equation*}
    V_\mu\cong\bigoplus_{\th\in P(X)}\Hom{Z_\la}{V_\mu}{\fr {g_0+t}}\otimes Z_\la
  \end{equation*}
The weight $\th$ will only appear in the $\la$-isotypic component if
$\la-\th\in R^+(X_0)$, the positive cone in the root lattice of $X_0$,
so the $\th$ weight space of $V_\mu$ can be decomposed as follows:
  \begin{equation*}
    \Hom{\C_{\th}}{V_{\mu}}{\fr t}=\bigoplus_{\xi\in
    R(X_0)}\Hom{\C_{\th}}{Z_{\th+\xi}}{\fr t}\otimes
    \Hom{Z_{\th+\xi}}{V_{\mu}}{\fr{g_0 + t}}
  \end{equation*}
  Taking dimensions, the result follows.
\end{proof} 

\begin{proof}[Proof of Theorem \ref{polyn-behav-weight}]
  In our notation, we can rewrite Lemma \ref{weight-decomp} as: 
\begin{equation}\label{eq:1}
    w^\Bw_\Bv( X(\bM))=\sum_{\Bz:\EV(X_0)\to\Z_{\geq 0}} b^{\Bw}_{\Bv-\Bz}
    (\bM) \cdot w^{\Bv_0+\Bz}_{\Bv_0}( X_0(\bM)). 
  \end{equation}

  While the set $\Znn^{\EV(X_0)}$ is infinite, the sum above is finite, since
  $b^{\Bw}_{\Bv-\Bz}(\bM)$ is non-zero for only finitely many
  $\Bz$.  

  Since
  $\Bv-\Bz\leq \Bv$, and $\Bv$ is of depth $s$, we can apply
  Corollary~\ref{vanishing-cor} to see that if $b^\Bw_{\Bv-\Bz}(X,
  X_0)\neq 0$, then $\Bv-\Bz$ is deep for $X_0^s$ of depth $\leq
  s$. 

  The set of weights deep on $X_0^s$ of depth $\leq s$ is independent of $m$ for large
  $m$, as is $b^{\Bw}_{\Bv-\Bz}(X, X_0)$. Thus, the set of $\Bz$ with
  $b^\Bw_{\Bv-\Bz}(X, X_0)\neq 0$ stabilizes as $\bM\to\infty$. 

  Thus, for large values of $\bM$, the sum in \eqref{eq:1} is over 
  a fixed set.  Its summands are polynomial since $b^{\Bw}_{\Bv-\Bz} 
  (\bM)$ stabilizes for all $y$, and by Theorem~\ref{benkart-et-al}, 
  the weight multiplicity $w^{\Bv_0+\Bz}_{\Bv_0}( X_0(\bM))$ is a
  polynomial of multidegree less than $\boldsymbol{\si}$. Thus,
  $w^{\Bw}_{\Bv} (\bM)$ is also a polynomial of multidegree $\leq
  \boldsymbol{\si}$ in the variables $m_1,\ldots, m_n$.
\end{proof}
\def\cprime{$'$}

\end{document}